\documentclass[oneside,draft,12pt]{amsart}
\usepackage{amssymb, amscd}
\usepackage[all]{xy}

 \newtheorem{theorem}{Theorem}[subsection]
 \newtheorem{corollary}[theorem]{Corollary}
 \newtheorem{lemma}[theorem]{Lemma}
 \newtheorem{proposition}[theorem]{Proposition}
 \theoremstyle{definition}
 \newtheorem{definition}[theorem]{Definition}
 \theoremstyle{definition}
 
 \theoremstyle{remark}
 \newtheorem{remark}[theorem]{Remark}
 \numberwithin{equation}{subsection}

\usepackage{latexsym}
\usepackage{amssymb}

\newcommand{\ben}{\begin{equation}}
\newcommand{\een}{\end{equation}}

\newcommand{\naturals}{\ensuremath{{\mathbb N}}}
\newcommand{\real}{\ensuremath{{\mathbb R}}}
\newcommand{\complex}{\ensuremath{{\mathbb C}}}

\newcommand{\rational}{\ensuremath{{\mathbb Q}}}

\newcommand{\GL}[1]{\ensuremath{{\mathrm {GL}_{ #1 }}}}
\newcommand{\GLC}[1]{\GL{#1}(\complex)}
\newcommand{\SL}[1]{\ensuremath{{\mathrm {SL}_{ #1 }}}}
\newcommand{\SLC}[1]{\SL{#1}(\complex)}

\newcommand{\OO}{\mathcal{O}}

\begin{document}
\title{Orbifold Cohomology of the Symmetric Product}

\author{Bernardo Uribe}
\address{Mathematics Department\\
         University of Wisconsin\\
         Madison, Wisconsin, 53706} 
\email{uribe@math.wisc.edu}

\begin{abstract}
Chen and Ruan's orbifold cohomology of the symmetric product
of a complex manifold is calculated. An
isomorphism of rings (up to a change of signs) $H_{orb}^*(X^n/S_n;\complex) \cong 
H^*(X^{[n]};\complex)$  between the orbifold cohomology
of the symmetric product of a smooth projective surface 
with trivial canonical class $X$ and the cohomology of its
Hilbert scheme $X^{[n]}$ is obtained, yielding a positive answer to
a conjecture of Ruan. 

\end{abstract}

\maketitle

\section{Introduction}

String theorists (see \cite{Dixon}) proposed an Euler characteristic for 
orbifolds that are global quotients by the action of a finite group;
 this number matched the
Euler characteristic of equivariant $K$-theory
(see \cite{AtiyahSegal}). For these global quotients,
 Zaslow \cite{Zaslow} worked out  an additive orbifold cohomology and later on, Chen and
Ruan \cite{ChenRuan} and Ruan \cite{Ruan} generalized this construction to 
a general orbifold. Motivated by the study
of quantum cohomology, they developed a new ring structure for the
cohomology of orbifolds. This ring structure is different to the ones
obtained  by other equivariant cohomology theories such as equivariant K-theory
 or Bredon Cohomology, and for the case of the
symmetric product it will be explicitly calculated in this paper.
For $X$ an algebraic surface, it is known that the Hilbert scheme $X^{[n]}$
of points of length
$n$ is a crepant resolution of $Symm_n(X)$; Ruan \cite{Ruan} conjectured that if 
$X^{[n]}$ had hyperk\"{a}hler structure then its cohomology and the orbifold
cohomology of $Symm_n(X)$ should be isomorphic as rings.  From
the calculation of the orbifold cohomology of $\complex^n/G$, with
$G$ a finite sobgroup of $GL(n,\complex)$, Ruan obtained
the orbifold cohomology of $Symm_n(\complex^2)$ which was matched
with the cohomology of $(\complex^2)^{[n]}$ obtained by Lehn
and Sorger \cite{LehnSorger2}. The same authors, in the case 
of a smooth projective surface with trivial canonical class $X$ calculated
the ring structure of the cohomology of $X^{[n]}$ \cite{LehnSorger}. Using 
the explicit calculation of the orbifold cohomology
of $Symm_n(X)$, when $X$ is a smooth projective surface with trivial canonical class,
an isomorphism (up to a change of sign) between the orbifold cohomology 
of $Symm_n(X)$ and the cohomology of $X^{[n]}$ is obtained.

The organization of the paper is as follows, in the first section the basic definitions
 of orbifold cohomology are summarized, 
and in the second section the cohomology ring structure of the symmetric product is
explained. After obtaining the results in this paper, I was notified 
that they had also been independently obtained by Fantechi and G\"{o}ttsche \cite{FantechiGottsche}.

 Last but not least, I would like to express my deepest gratitude
to E. Lupercio, M. Poddar, A. Adem and Y. Ruan who shared with me insightful ideas
through informal meetings and especially to the latter two who introduced
me to the subject.

\section{Preliminaries}

\subsection{Orbifold Cohomology}

For a full treatment of orbifolds the papers of Chen and Ruan \cite{ChenRuan}
and of Ruan \cite{Ruan} are recommended. The definitions, notations and results 
of those papers will be used in what follows.
Let $Y=X/G$ be a global quotient by a finite group, thus in particular an orbifold. $T_k$
will be the set of conjugacy classes of $k$-tuples ${\bf g}=(g_1,\dots,g_k)$ 
of elements in $G$.

The twisted sectors are the sets
$$Y_{({\bf g})} =X^{{\bf g}}/C({\bf g})$$
where $X^{{\bf g}}=X^{g_1} \cap \cdots \cap X^{g_k}$  and
$C({\bf g})=C(g_1)\cap \cdots \cap C(g_k)$, where $C(g_i)$ is the centralizer of $g_i$ in $G$.
The multisectors $\widetilde{\Sigma_k Y}$ are the disjoint union of the twisted sectors, i.e.
$$\widetilde{\Sigma_k Y}= \bigsqcup_{({\bf g}) \in T_k} Y_{({\bf g})}$$

Let's now consider the natural maps between multi-sectors; the
evaluation maps $e_{i_1, \dots, i_l} : \widetilde{\Sigma_k Y} \to
\widetilde{\Sigma_l Y}$ defined by $e_{i_1, \dots,
i_l}(x,(g_1,\dots, g_k)) \mapsto (x,(g_{i_1}, \dots, g_{i_l}))$
and the involutions $I:\widetilde{\Sigma_k Y} \to
\widetilde{\Sigma_k Y}$ defined by $I(x,({\bf g})) \mapsto
(x,({\bf g^{-1}}))$ where ${\bf
g^{_1}}=(g_1^{-1},\dots,g_k^{-1})$.

To define the orbifold cohomology group we need to add a shifting
to  the cohomology of the twisted sectors, and for that we are
going to assume that the orbifold $Y$ is almost complex with
complex structure $J$; recall that $J$ will be a smooth section
of $End(TY)$ such that $J^2 =-Id$.

For $p \in Y$ the almost complex structure gives rise to a
faithful representation $\rho_p : G \to \GLC{n}$ ($n =
dim_{\complex} Y$) that could be diagonalized as
$$diag \left(e^{2 \pi \frac{m_{1,g}}{m_g}}, \dots ,
 e^{2 \pi \frac{m_{n,g}}{m_g}} \right)$$
where $m_g$ is the order of $g$ in $G$ and $0 \leq m_{j,g} <
m_g$. We define a function $\iota : \widetilde{ \Sigma Y} \to
\rational$ by
$$\iota(p,(g)) = \sum_{j=1}^n \frac{m_{j,g}}{m_g}$$
It is easy to see that it is locally constant, hence we call it
$\iota_{(g)}$, the degree shifting number in each sector $Y_{(g)}$;
 it is an integer if and only if $\rho_p(g) \in
\SLC{n}$ and
$$\iota_{(g)} + \iota_{(g^{-1})} = rank(\rho_p(g) - I)$$
which is the complex codimension $dim_{\complex} Y -
dim_{\complex} Y_{(g)}$.

\begin{definition}
The orbifold cohomology
groups are defined as
$$H_{orb}^d(Y; \complex) = \bigoplus_{(g) \in T_1} H^{d- 2 \iota_{(g)}}
 (Y_{(g)}; \complex)$$
\end{definition}

\subsection{Poincar\'{e} Duality}

Let's assume the orbifold $Y$ is closed and recall the map $I:
Y_{(g)} \to Y_{(g^{-1})}$ defined by $(p,(g)) \mapsto
(p,(g^{-1}))$ with $I^2=Id$.

Poincar\'{e} Duality for the orbifold cohomology is as follows:

For any $0 \leq d \leq 2n$, the pairing
$$<>_{orb}: H^d_{orb}(Y;\complex) \otimes H^{2n-d}_{orb}(Y;\complex) \to
\complex$$ defined by the direct sum of
$$<>^{(g)}_{orb}: H^{d-2 \iota_{(g)}}(Y_{(g)};\complex) \otimes
H^{2n- d-2 \iota_{(g^{-1})}}(Y_{(g^{-1})};\complex) \to
\complex$$ where $$<\alpha, \beta>^{(g)}_{orb} :=
\int_{Y_{(g)}}\alpha \wedge I^*(\beta)$$ is nondegenerate and
$\alpha \in H^{d-2 \iota_{(g)}}(Y_{(g)};\complex)$, $\beta \in
H^{2n- d-2 \iota_{(g^{-1})}}(Y_{(g^{-1})};\complex)$.

When restricted to the non-twisted sector this is the
ordinary Poincar\'{e} pairing.

\subsection{Orbifold Cup Product }

The orbifold cup product relies on the construction of a
obstruction bundle over the twisted sector $Y_{({\bf g})}$ where
$({\bf g}) =(g_1, g_2, g_3) \in T_3$ is the conjugacy class of
the triple $(g_1, g_2, g_3)$ with $g_1g_2g_3 =1$, and $T_3^0$ is
the set of those conjugacy classes.

Let $e: Y_{({\bf g})} \to Y$ be the evaluation map and $e^*TY$ the
pullback tangent bundle over $Y_{({\bf g})}$. For $y \in Y_{({\bf
g})}$ its local group in $Y$ is $\Gamma'$; clearly it contains the
elements $g_1, g_2, g_3$ with the relations $g_1g_2g_3=1$ and
$g_i^{k_i}$, where $k_i$ is the order of $g_i$. Let $\Gamma$ be the
subgroup of $\Gamma'$ generated by these three elements $g_1, g_2,
g_3$; then $\Gamma$ acts on $e^*TY$ while fixing $Y_{({\bf g})}$.

For the orbifold Riemann sphere with three orbifold points $(S^2,
(x_1,x_2,x_3),(k_1,k_2,k_3))$ there exists a closed Riemann
surface $\Sigma$ such that $\Gamma$ acts on it holomorphically,
$\Sigma /\Gamma = S^2$ and $(\Sigma,\Gamma,\pi)$ is a uniformizing system
for $S^2$ . The group $\Gamma$ acts on both
$H^1(\Sigma)$ and $e^*TY$ where we consider $H^1(\Sigma)$ as a
trivial bundle over $Y_{({\bf g})}$.

The obstruction bundle $E_{({\bf g})}$ we require is the invariant
part of $H^1(\Sigma) \otimes e^*TY$ under the action of $\Gamma$, i.e.
$E_{({\bf g})} = \left( H^1(\Sigma) \otimes e^*TY \right) ^\Gamma$.
Let $c(E_{({\bf g})})$ be the Euler class of $E_{({\bf g})}$
(which up to an exact form is independent of the connection) and
recall the evaluation maps $e_i : Y_{({\bf g})} \to Y_{(g_i)}$.

\begin{definition}
For $\alpha, \beta,\gamma \in H^*_{orb}(Y;\complex)$ a three-point
function is defined
$$<\alpha,\beta,\gamma>_{orb} = \sum_{({\bf g}) \in T_3^0}
\int_{Y_{({\bf g})}} e_1^*\alpha \cdot e_2^*\beta \cdot
e_3^*\gamma \cdot c(E_{({\bf g})})$$ and let the orbifold cup
product be defined by the relation
$$<\alpha \cup_{orb} \beta, \gamma >_{orb} =
<\alpha,\beta,\gamma>_{orb}$$
\end{definition}

\begin{lemma}
For $\alpha \in H^*(Y_{(g_1)}; \complex)$ and $\beta \in
H^*(Y_{(g_2)}; \complex)$ the cup product $\alpha \cup_{orb}
\beta$ can be decomposed as a sum of its components in
$H^*_{orb}(Y; \complex)= \bigoplus_{(g)} H^*(Y_{(g)};\complex)$. So
we have
$$\alpha \cup_{orb} \beta = \sum_{ \begin{scriptsize}
\begin{array}{c}
(h_1,h_2, (h_1h_2)^{-1}) \in T_3^0 \\
h_i \in (g_i)
\end{array} \end{scriptsize}}
p_{*} \left( e_1^*\alpha \cdot  e_2^*\beta \cdot c(E_{({\bf
h})}) \right)$$ where $p_{*}$ is the push-out of the evaluation
map $p: Y_{(h_1,h_2)} \to Y_{(h_1h_2)}$.
\end{lemma}

\begin{remark} \label{remarkthreepointfunction}
In the definition of the three point function for the orbifold
cup product an abuse of notation is being made. For $\alpha \in
H^*(Y_{(g_1)};\complex)$ we need to take into account the
different elements conjugated to $g_1$, for $h_1 \in (g_1)$,
$\phi:H^*(Y_{(g_1)};\complex) \stackrel{\cong}{\to}
H^*(Y_{(h_1)};\complex)$ and let's denote by
$\alpha_{h_1}:=\phi(\alpha)$. Then for $\beta \in
H^*(Y_{(g_2)};\complex)$ and $\gamma \in H^*(Y_{(g_3)};\complex)$
$$<\alpha,\beta,\gamma>_{orb} = \sum_{
\begin{scriptsize}
\begin{array}{c}
{({\bf h}) \in T_3^0} \\
h_i \in (g_i) \end{array} \end{scriptsize}}
 \int_{Y_{(h_1,h_2)}}
e_1^*\alpha_{h_1} \cdot e_2^*\beta_{h_2} \cdot e_3^*\gamma_{h_3}
\cdot c(E_{({\bf h})})$$
\end{remark}

\bigskip

 The previous description becomes
very complicated when calculations are tried out. In what follows,
an equivalent description for the orbifold cohomology that
simplifies calculations will be explained. The idea comes from the paper of Lehn
and Sorger \cite{LehnSorger} and similar notation will be used.

Let $Y=X/G$ be an orbifold with $X$ a compact
complex manifold. As before $X^g$ will denote the fixed
point set of the action of $g$ on $X$.

The cohomology classes will be labeled by elements in $G$ and let
the total ring $A(X,G)$ be
$$A(X,G) := \bigoplus_{g \in G} H^*(X^g;\complex) \times \{g\}$$
Its group structure is the natural one and the ring structure
that will be defined later will give us the orbifold cup product.
The grading is the one  in the orbifold cohomology, i.e.
$$A^d(X,G) = \bigoplus_{g \in G} H^{d-2\iota_{(g)}}(X^g;\complex) \times \{g\}$$

For $h \in G$ there is a natural map $h: X^g \to X^{hgh^{-1}}$
which can be extended to an action in $A(X,G)$ inducing an
isomorphism \begin{eqnarray*} h :H^*(X^g;\complex) \times \{g\} &
\to & H^*(X^{hgh^{-1}};\complex) \times \{hgh^{-1}\}\\
(\alpha,g) & \mapsto & ((h^{-1})^*\alpha, hgh^{-1})
\end{eqnarray*}

The invariant part under the action of $G$ is isomorphic as a group to the
orbifold cohomology,

\begin{lemma} \label{isomorphismofgroups}
$$A(X,G)^G \cong \bigoplus_{(g)} H^*(X^g;\complex)^{C(g)} \cong
H^*_{orb}(X/G;\complex)$$
\end{lemma}

Now some notation needs to be introduced in order to define the
ring structure. This notation comes from \cite{LehnSorger} where
a more detailed study of Frobenius algebras is done.

\begin{definition}
For $X^{\langle h_1,h_2 \rangle}$, the fixed point set of $\langle h_1,h_2 \rangle$,
 let $$f^{h_i,\langle h_1,h_2 \rangle} :
H^*(X^{h_i};\complex) \to H^*(X^{\langle h_1, h_2 \rangle};\complex)$$
$$f_{\langle h_1,h_2 \rangle, h_i} :
H^*(X^{\langle h_1, h_2 \rangle};\complex) \to H^*(X^{h_i};\complex)$$ be the
pull-back and the push-forward respectively of the diagonal inclusion
map $X^{\langle h_1, h_2 \rangle} \hookrightarrow X^{h_i}$ where $i=1,2,3$ and
$h_3=(h_1h_2)$.
\end{definition}

We need to make use of the obstruction bundle over $X^{\langle h_1, h_2 \rangle}$;
 as $Y_{(h_1,h_2)} = X^{\langle h_1, h_2 \rangle}/C(h_1,h_2)$ and taking the
 projection map $\pi :X^{\langle h_1, h_2 \rangle} \to X^{\langle h_1, h_2 \rangle}/C(h_1,h_2)$ we
 will consider the Euler class of $\pi^*(E_{({\bf h})})$.

 \begin{definition}
Let the product $A(X,G) \otimes A(X,G) \stackrel{\cdot}{\to}
A(X,G)$ be defined by
$$(\alpha,h_1) \cdot (\beta,h_2) := f_{\langle h_1,h_2 \rangle, h_1h_2} \left(
f^{h_1,\langle h_1,h_2 \rangle} (\alpha) \wedge f^{h_2,\langle h_1,h_2
\rangle}(\beta)  \wedge \pi^*c(E_{({\bf h})}) \right)$$
whose three point function is
\begin{eqnarray*}
\lefteqn{<(\alpha,h_1),(\beta,h_2),(\gamma,(h_1h_2)^{-1})>} & &\\
& := & \int_{X^{\langle h_1, h_2 \rangle}}
 f^{h_1,\langle h_1,h_2 \rangle} (\alpha) \wedge f^{h_2,\langle h_1,h_2
\rangle}(\beta)  \wedge f^{(h_1h_2)^{-1}, \langle h_1,h_2
\rangle}(\gamma) \wedge \pi^*c(E_{({\bf h})})
\end{eqnarray*}
 \end{definition}

\begin{lemma}
The product $A(X,G) \otimes A(X,G) \stackrel{\cdot}{\to} A(X,G)$
previously defined is $G$ equivariant.
\end{lemma}

This product induces a ring structure on the invariant group
$A(X,G)^G$ which will match with the orbifold cup product. Thus $A(X,G)^G$
will inherit the properties of the orbifold cup product.

\begin{proposition}
The rings  $H^*_{orb}(X/G;\complex)$ and $A(X,G)^G$ are
isomorphic.
\end{proposition}

\begin{proof}
From lemma
\ref{isomorphismofgroups} it is known that they are isomorphic as groups,
but we need to have an explicit map
$$\varphi:H^*_{orb}(X/G,\complex) \to A(X,G)^G.$$
 For $\alpha \in
H^*(Y_{(g_1)};\complex)$ let $\alpha' \in
H^*(X^{g_1};\complex)^{C(g_1)}$ be the pullback of $\alpha$ under
the projection map $X^{g_1} \to Y_{(g_1)}=X^{g_1}/C(g_1)$. In the
same way as in remark \ref{remarkthreepointfunction} let
$(\alpha')_{h_1} := \phi'(\alpha')$ where $\phi':
H^*(X^{g_1};\complex) \stackrel{\cong}{\to}
H^*(X^{h_1};\complex)$ and $h_1 \in (g_1)$. As $(\alpha')_{h_1} =
(\alpha_{h_1})'$ there will be no confusion in denoting it by
$\alpha'_{h_1}$.

The isomorphism of groups is clearly given by
$$\varphi(\alpha):= \sum_{h_1 \in (g_1)} (\alpha'_{h_1},h_1)$$
and what is left to prove is that the triple functions give the
same result.

For $\alpha \in H^*(Y_{(g_1)},\complex)$, $\beta \in
H^*(Y_{(g_2)},\complex)$ and $\gamma \in H^*(Y_{(g_3)},\complex)$
the triple function $<>_G$ induced by $<>$ in $A(X,G)^G$ is
\begin{eqnarray*}
\lefteqn{<\varphi(\alpha),\varphi(\beta),\varphi(\gamma)>_G} &
& \\
& = & \frac{1}{|G|} \sum_{
\begin{scriptsize}
\begin{array}{c}
h_i \in (g_i) \\
h_3 = (h_1h_2)^{-1}
\end{array}
\end{scriptsize}}
\int_{X^{\langle h_1, h_2 \rangle}} f^{h_1,\langle h_1,h_2 \rangle}(\alpha'_{h_1})
\cdot f^{h_2,\langle h_1,h_2 \rangle}(\beta'_{h_2}) \cdot
f^{h_3,\langle h_1,h_2 \rangle}(\gamma'_{h_3}) \cdot
\pi^*c(E_{({\bf h})})
\end{eqnarray*}
The group $G$ acts on the integral via conjugation, and it is clear
that its value is invariant under this action. The action of $G$
by conjugation on the pairs $[h_1,h_2]$ has as stabilizer the
group $C(h_1,h_2)$, thus we can exchange the set of pairs
$[h_1,h_2]$ with $h_i \in (g_i)$ by the set of conjugacy classes
of pairs $(h_1,h_2)$ multiplying  by $|G|/|C(h_1,h_2)|$, the
multiplicity of the integral ; hence

\begin{eqnarray*}
\lefteqn{<\varphi(\alpha),\varphi(\beta),\varphi(\gamma)>_G} &
& \\
& = & \sum_{
\begin{scriptsize}
\begin{array}{c}
({\bf h}) \in T_3^0 \\
 h_i \in (g_i)
\end{array}
\end{scriptsize}}
\frac{1}{|C(h_1,h_2)|}\int_{X^{\langle h_1, h_2 \rangle}} f^{h_1,\langle h_1,h_2
\rangle}(\alpha'_{h_1}) \cdot f^{h_2,\langle h_1,h_2
\rangle}(\beta'_{h_2}) \cdot f^{h_3,\langle h_1,h_2
\rangle}(\gamma'_{h_3}) \cdot \pi^*c(E_{({\bf h})})
\end{eqnarray*}

Changing the set of integration and using that $$\int_{Y_{(g,h)}}
\omega = \frac{1}{|C(g,h)|} \int_{X^{\langle g,h \rangle}} \pi^* \omega$$ we get
the desired equality

\begin{eqnarray*}
\lefteqn{<\varphi(\alpha),\varphi(\beta),\varphi(\gamma)>_G} &
& \\
& = & \sum_{
\begin{scriptsize}
\begin{array}{c}
({\bf h}) \in T_3^0 \\
 h_i \in (g_i)
\end{array}
\end{scriptsize}}
\int_{Y_{(h_1,h_2)}} e_1^*(\alpha_{h_1}) \cdot e_2^*(\beta_{h_2})
\cdot e_3^*(\gamma_{h_3}) \cdot c(E_{({\bf h})}) \\
 & = & <\alpha,\beta,\gamma>_{orb}
\end{eqnarray*} where one recalls that
$\alpha'_{h_1}=\pi^*\alpha_{h_1}$.

 As the triple functions give
the same result and the grading matches in both descriptions,
 the ring isomorphism follows. All the properties of the
orbifold cup product proven in \cite{ChenRuan} apply to
$A(X,G)^G$; this will be the description of orbifold cohomology
that will be useful in the calculations that follow.

\end{proof}

\section{Orbifold Cohomology of the Symmetric Product}

In this section the previous description of the orbifold cup
product of global quotients will be applied to the symmetric
product. $X$ will be an even dimensional complex manifold 
 $dim_{\complex} X =2N$, and the orbifold in
mind will be $X^n/S_n$ where the action of the symmetric group
$S_n$ on $X^n$ is the natural one.

Some notation needs to be introduced.

{\bf Notation}: For $\sigma, \rho \in S_n$, let $\Gamma \subset
[n]:=\{1,2,\dots,n\}$ be a set stable under the action of
$\sigma$; we will denote by $\OO(\sigma;\Gamma)$ the set of orbits
induced by the action of $\sigma$ in $\Gamma$. If $\Gamma$ is
$\sigma$-stable and $\rho$-stable, $\OO(\sigma,\rho;\Gamma)$ will
be the set of orbits induced by $\langle \sigma,\rho \rangle$.
When the set $\Gamma$ is dropped from the expression, the
set $\OO(\sigma, [n])$ will be  denoted $\OO(\sigma)$.

$|\sigma|$  will denote  the minimum number $m$ of transpositions
$\tau_1, \dots ,\tau_m$ such that $\sigma=\tau_1\dots\tau_m$; hence
$$|\sigma|+|\OO(\sigma)|=n$$ 

The set $X^n_\sigma$ will denote the fixed point set under the action of
$\sigma$ on $X^n$. Superscripts on $X$ will count the number of copies of itself
on the cartesian product, and subscripts will be elements of the group and
will determine fixed point sets.

 For $A$ a graded ring and $m$ an
integer, we will denote by $A[m]$ the ring $A$ whose grading is
being shifted $m$ units to the left; in other words
$$A[m]^{i}:=A^{i-m}$$

 The orbifold ring structure is easy to understand from $A(X^n,S_n)^{S_n}$ following
the construction given previously. The only things left to study
are the obstruction bundles and some splitings that occur using
the orbits of the actions in $[n]$ of the elements of $S_n$.

\begin{lemma}
The orbifold shifting number of $\sigma \in S_n$ is
$$\iota_{(\sigma)}=\sum_j n_j \frac{(j-1)}{2}N = \frac{N}{2}|\sigma|$$
so that $2\iota_{(\sigma)} =N|\sigma|$.
\end{lemma}
\begin{proof}
Once the action of $\sigma$ in $X^n$ is diagonalized every
$j$-cycle gives all the $j$-roots of unity; the sum of these roots
is $\frac{(j-1)}{2}N$.
\end{proof}
As $|\sigma| + |\OO(\sigma)|=n$ the shifting is taking into
account that $X^n_{\sigma}$ is isomorphic to $X^{|\OO(\sigma)|}$.
Via the Kunneth isomorphism we have the following set of
identities (everywhere the coefficient system for cohomology will
be $\complex$, so it will be dropped out of the notation)
$$H^*(X^n_{\sigma}) \cong H^*(X)^{\otimes |\OO(\sigma)|}$$
then
$$H^*(X)[N]^{\otimes |\OO(\sigma)|} \cong
H^*(X^n_{\sigma})[|\OO(\sigma)|N] \cong
H^{*-2\iota_{(\sigma)}}(X^n_{\sigma}) [nN]$$

\medskip
 Using the notation of Lehn and Sorger \cite[Section 2]{LehnSorger}
 where
$$H^*(X)[N]\{S_n\} = \bigoplus_{\sigma \in S_n}
H^*(X)[N]^{\otimes |\OO(\sigma)|} \cdot \sigma$$ and the previous
isomorphisms we get

\begin{proposition}
$H^*(X)[N]\{S_n\}$ and $A(X^n,S_n)[nN]$ are isomorphic as graded
vector spaces.
\end{proposition}

The ring structure of $H^*(X)[N]\{S_n\}$ is defined by the
product (see \cite{LehnSorger})

$$m_{\pi,\rho}: H^*(X)[N]^{\otimes|\OO(\pi)|} \otimes H^*(X)[N]^{\otimes|\OO(\rho)|}
\to H^*(X)[N]^{\otimes|\OO(\pi \rho)|}$$
$$m_{\pi,\rho}(\alpha \otimes \beta) = f_{\langle \pi ,\rho \rangle, \pi\rho} \left(
f^{\pi,\langle \pi,\rho \rangle} (\alpha) \wedge f^{\rho,\langle
\pi,\rho \rangle}(\beta)  \wedge e^{g(\pi,\rho)} \right)$$ where
$g(\pi,\rho): \OO(\pi,\rho) \to \naturals$ is the graph defect
function; for $\Gamma \in \OO(\pi,\rho)$
\begin{eqnarray} \label{graphdefectfunction}
g(\pi,\rho)(\Gamma)= \frac{1}{2} \left(
|\Gamma| +2 -|\OO(\pi,\Gamma)| -|\OO(\rho,\Gamma)| -
|\OO(\pi\rho,\Gamma) \right)
\end{eqnarray}
 and
$$e^{g(\pi,\rho)}:= \prod_{\Gamma \in \OO(\pi,\rho)}
e(X)^{g(\pi,\rho)(\Gamma)}$$ where $e(X)$ is the Euler class of $X$.

For simplicity an abuse of notation is being made; the functions 
$f^{\pi,\langle \pi,\rho \rangle}$ are used to define either one
of the following morphisms 

\begin{displaymath}
\begin{array}{ccccc}
H^*(X^n_{\pi}) & & \stackrel{f^{\pi,\langle \pi,\rho \rangle}}{\longrightarrow}
& & H^*(X^n_{\pi,\rho})\\
\cong \updownarrow & & & & \cong \updownarrow\\
H^*(X)^{\otimes|\OO(\pi)|} & & \stackrel{f^{\pi,\langle \pi,\rho \rangle}}{\longrightarrow}
& & H^*(X)^{\otimes|\OO(\pi,\rho)|}\\
\end{array}
\end{displaymath}

It is clear that the product in  $A(X^n,S_n)$ and in
$H^*(X)[N]\{S_n\}$are defined almost identically. The difference is
on the last term, the Euler class of the obstruction bundle on
one side and $e^{g(\pi,\rho)}$ on the other. We will see that
these two terms represent the same class.

\subsection{The obstruction bundle}
For $h_1,h_2 \in S_n$ the obstruction bundle $E_{({\bf h})}$ over
$Y_{(h_1,h_2)}$ is defined by
$$E_{({\bf h})}= \left( H^1(\Sigma) \otimes e^*TY \right)^G$$
where $G= \langle h_1,h_2 \rangle$, $Y=X^n/S_n$ and $\Sigma$ is an
orbifold Riemann surface provided with a $G$ action such that
$\Sigma / G =(S^2, (x_1,x_2,x_3),(k_1,k_2,k_3))$ is an orbifold
sphere with three marked points.

Because $H^1(\Sigma)$ is a trivial bundle, the pullback of
$E_{({\bf h})}$ under $\pi: X^n_{h_1,h_2} \to Y_{(h_1.h_2)}$ is
$$E_{h_1,h_2}:=\pi^*E_{({\bf h})} = \left( H^1(\Sigma) \otimes \Delta^*TX^n \right)^G$$
where $\Delta: X^n_{h_1,h_2} \hookrightarrow X^n$ is the inclusion (if $\rho : X^n \to Y$
is the quotient map, then $\rho \circ \Delta = e \circ \pi$).

Without loss of generality we can assume that $|\OO(h_1,h_2)|=k$
and $n_1+ \cdots + n_k =n$ a partition of  such that
$$\Gamma_i =\{n_1 + \cdots + n_{i-1}+1, \dots , n_1 + \cdots +n_i
\}$$ and $\{ \Gamma_1, \Gamma_2, \dots , \Gamma_k \} =
\OO(h_1,h_2)$. We will concentrate on each of the $\Gamma_i$
because we will see that the obstruction bundle $E_{h_1,h_2}$ can
be seen as the product of $k$ bundles over $X$ (i.e. $E_{h_1,h_2}
= \prod_i E_{h_1,h_2}^i$).

The following commutative diagrams
$$
\xy\xygraph{
!M{ X^n  & X^n_{h_1} \:[l]  \\
X^n_{h_2} \:[u]  & X^n_{h_1,h_2}\:[l] \:[u] } }\endxy
\ \ \ \ \ \ \ \ \ 
\xy\xygraph{
!M{ [n] \:[r] \:[d] & \OO(h_1) \:[d]  \\
\OO(h_2) \:[r]  & \OO(h_1,h_2) } }\endxy
$$

where arrows in the first one are inclusions of sets and the second are
inclusion of orbits, induce the commutative diagram of diagonal inclusions
$$
\xy\xygraph{
!M{ X^n  & X^{|\OO(h_1)|} \:[l]  \\
X^{|\OO(h_2)|} \:[u]  & X^{|\OO(h_1,h_2)|} \:[l] \:[u] } }\endxy
$$
where every arrow is the product of the arrows of the following diagram
for $i=1,2,\dots,k$
$$
\xy\xygraph{
!M{ X^{n_i} \:[r] \:[d] & X^{|\OO(h_1; \Gamma_i)|} \:[d]  \\
X^{|\OO(h_2;\Gamma_i)|} \:[r]  & X^{|\OO(h_1,h_2;\Gamma_i)|}=X } }\endxy
$$

\begin{lemma}
For $\Delta_i : X \to X^{n_i}$ $i=1,\dots,k$  the diagonal inclusions,
the bundles $\Delta_i^*TX^{n_i}$ become $G$ bundles via the
restriction of the action of $G$ into the orbit $\Gamma_i$ and
$$\Delta^*TX^n \cong \Delta_1^*TX^{n_1} \times \cdots \times
\Delta_k^*TX^{n_k}$$ as $G$ vector bundles.
\end{lemma}

\begin{proof} This comes from the fact that the orbits $\Gamma_i$
are $G$ stable, hence $G$ induces an action on each $X^{n_i}$.
\end{proof}

\begin{corollary}
The obstruction bundle splits as
$$E_{h_1,h_2} = \prod_{i=1}^k \left( H^1(\Sigma) \otimes
\Delta_i^*TX^{n_i} \right) ^G$$
\end{corollary}

We can simplify  the previous expression a bit further. Let $G_i$ be
 the subgroup of $S_{n_i}$ obtained from $G$ when its action is restricted to
the elements in $\Gamma_i$; then we have a surjective homomorphism
$$\lambda_i : G \to G_i$$ where the action of $G$ into
$\Delta_i^*TX^{n_i}$ factors through $G_i$. So we have

\begin{lemma}
$\left( H^1(\Sigma) \otimes \Delta_i^*TX^{n_i} \right) ^G \cong
\left( H^1(\Sigma)^{ker(\lambda_i)} \otimes \Delta_i^*TX^{n_i}
\right) ^{G_i}$
\end{lemma}

Now let $\Sigma_i:= \Sigma / \ker(\lambda_i)$, it is an orbifold
Riemann surface with a $G_i$ action so that  $\Sigma_i /G_i$
becomes an orbifold sphere with three marked points (the markings
are with respect to the generators of $G_i$:
$h_1\upharpoonright_{\Gamma_i}, h_2\upharpoonright_{\Gamma_i}$ and $
h_1h_2\upharpoonright_{\Gamma_i}$). So, in the same way as in the
definition of the obstruction bundle $E_{({\bf h})}$ we get that
$$E_{h_1,h_2}^i:=\left( H^1(\Sigma_i) \otimes \Delta_i^*TX^{n_i}
\right) ^{G_i}$$ hence

\begin{proposition}
The obstruction bundle splits as
$$E_{h_1,h_2}=\prod_{i=1}^k E_{h_1,h_2}^i$$
\end{proposition}

As the action of $G_i$ in $\Delta_i^*TX^{n_i}$ is independent on the
structure of $X$( moreover, it depends only in the coordinates),
hence

$$\Delta_i^*TX^{n_i} \cong TX \otimes \complex^{n_i}$$
as $G_i$-vector bundles, where $TX$ is the tangent bundle over $X$
and $G_i \subset S_{n_i}$ acts on $\complex^{n_i}$ in the natural
way. Then

\begin{lemma}
$$E_{h_1,h_2}^i \cong TX \otimes (H^1(\Sigma) \otimes
\complex^{n_i})^{G_i}$$
\end{lemma}

Defining $r(h_1,h_2)(i):=dim_\complex (H^1(\Sigma) \otimes
\complex^{n_i})^{G_i}$ it follows that the Euler class of
$E_{h_1,h_2}^i$ equals the Euler class of $X$ to some exponent

\begin{corollary}
$c(E_{h_1,h_2}^i) = e(X)^{r(h_1,h_2)(i)}$
\end{corollary}
in other words,
$$c(E_{h_1,h_2}^i) = \left\{ \begin{array}{cc}
1 & \mbox{if }{r(h_1,h_2)(i)}=0 \\
e(X) & \mbox{if } {r(h_1,h_2)(i)}=1\\
0 & \mbox{if }{r(h_1,h_2)(i)} \geq 2
 \end{array} \right. $$

Making it clear that the obstruction looks like
$$c(E_{h_1,h_2}) = \prod_{i=1}^k e(X)^{r(h_1,h_2)(i)}$$

\begin{proposition}
$H^*(X)[N]\{S_n\}$ and $A(X^n,S_n)[nN]$ are isomorphic as rings.
\end{proposition}
\begin{proof}
The only thing left to prove is that the graph defect function
$g(h_1,h_2)$ matches the function $r(h_1,h_2)$ just defined.
 Working on the orbifold cup product, using the previous commutative
 diagrams (especially \ref{comm.diagram.of.each.orbit}) and the splitings we can
 restrict ourselves to each orbit $\Gamma_i$ of $\OO(h_1,h_2)$,
 hence for
 $$\alpha \in
 H^{p-2\iota_{(h_1\upharpoonright_{\Gamma_i})}}(X^{|\OO(h_1;\Gamma_i)|})$$
$$\beta \in
 H^{q-2\iota_{(h_2\upharpoonright_{\Gamma_i})}}(X^{|\OO(h_2;\Gamma_i)|})$$
 $$\alpha \in
 H^{r-2\iota_{((h_1h_2)^{-1}\upharpoonright_{\Gamma_i})}}(X^{|\OO((h_1h_2)^{-1};\Gamma_i)|})$$
 with $p+q+r=2n_iN$, and $|\Gamma_i|=n_i$
 \begin{eqnarray*}
deg \,\alpha + deg \,\beta + deg \,\gamma & = & p+q+r
-2(\iota_{(h_1\upharpoonright_{\Gamma_i})}
+\iota_{(h_2\upharpoonright_{\Gamma_i})}+\iota_{((h_1h_2)^{-1}\upharpoonright_{\Gamma_i})})\\
& = & (2n_i -|h_1\upharpoonright_{\Gamma_i}| -|h_2\upharpoonright_{\Gamma_i}|
 -|h_1h_2\upharpoonright_{\Gamma_i}|)N \\
& = & (|\OO(h_1;\Gamma_i)| +|\OO(h_2;\Gamma_i)|
+|\OO(h_1h_2;\Gamma_i)| - n_i)N
 \end{eqnarray*}
so we get that
\begin{eqnarray*}
dim_{\real} E_{h_1,h_2}^i & = & dim_{\real} X - deg \,\alpha -deg
\,\beta - deg \,\gamma\\
& = & \frac{1}{2}(2 + n_i  - |\OO(h_1;\Gamma_i)|
-|\OO(h_2;\Gamma_i)| -|\OO(h_1h_2;\Gamma_i)|) 2N
\end{eqnarray*}

which implies that
$$r(h_1,h_2)(i)= \frac{1}{2}(2 + n_i  - |\OO(h_1;\Gamma_i)|
-|\OO(h_2;\Gamma_i)| -|\OO(h_1h_2;\Gamma_i)|)$$ matching
precisely the definition of the graph defect function
$g(h_1,h_2)$ (see formula \ref{graphdefectfunction}).

\end{proof}

Using the notation of \cite{LehnSorger}
$$H^*(X)[N]^{[n]}:=(H^*(X)[N]\{S_n\})^{S_n}$$
we can conclude

\begin{theorem} \label{thmorbcohomology}
For $X$ a compact complex even dimensional manifold
($dim_\complex (X)=N$) we have that
$$H^*_{orb}(X^n/S_n;\complex)[nN] \cong H^*(X)[N]^{[n]}$$
\end{theorem}
\begin{proof}
The rings $H^*(X)[N]\{S_n\}$ and $A(X^n,S_n)[nN]$ are isomorphic
and they have a compatible $S_n$ equivariant action, hence their
invariant rings are also isomorphic, so
$$H^*_{orb}(X^n/S_n;\complex)[nN] \cong A(X^n,S_n)[nN]^{S_n} \cong
(H^*(X)[N]\{S_n\})^{S_n}$$
\end{proof}
\subsection{Hilbert Schemes}

Now we are ready to prove a conjecture posed by Ruan \cite[Conj.
6]{Ruan} about algebraic surfaces. Let $X$ be a smooth projective
surface over the complex numbers and $X^{[n]}$ its
$n$-th Hilbert scheme of points of length $n$; $X^{[n]}$ is again projective and smooth of
dimension $2n$. Lehn and Sorger proved
\begin{theorem} \cite[Thm. 3.2]{LehnSorger} \label{thmlehnsorger}
Let $X$ be a smooth projective surface with trivial canonical
divisor. Then there is a canonical isomorphism of graded rings
$$(H^*(X;\complex)[2])^{[n]} \cong H^*(X^{[n]};\complex)[2n]$$
\end{theorem}

\begin{remark}
This isomorphism is obtained after changing the sign on the integral over
the fundamental homology class of $X$. In the notation of \cite{LehnSorger}
$T(a) :=- \int_{[X]}a$, for $a \in H^*(X;\complex)$.
\end{remark}
Using the results of the previous section we obtain a positive
answer to the conjecture

\begin{theorem}
Let $X$ be a smooth projective surface with trivial canonical
divisor. Then there is a canonical isomorphism (up to a change of sign)
of graded rings  between the orbifold cohomology of
$X^n/S_n$ and the cohomology of the $n$-th Hilbert scheme of $X$
$$H^*_{orb}(X^n/S_n;\complex) \cong H^*(X^{[n]};\complex)$$
\end{theorem}
\begin{proof}
By theorems \ref{thmlehnsorger} and \ref{thmorbcohomology} 
$$H^*_{orb}(X^n/S_n;\complex)[2n] \cong (H^*(X;\complex)[2])^{[n]}
\cong H^*(X^{[n]};\complex)[2n]$$

\end{proof}

\bibliographystyle{amsplain}

\providecommand{\bysame}{\leavevmode\hbox
to3em{\hrulefill}\thinspace}

\end{document}